\documentclass[letterpaper,12pt]{article}

\usepackage{graphicx}
\usepackage{amssymb,amsmath}
\usepackage{subfigure}
\usepackage{amsthm}
\usepackage{hyperref}
\usepackage{color}
\usepackage{multirow,txfonts,verbatim}
\usepackage{mathrsfs}
\usepackage{booktabs}
\usepackage{bm}
\usepackage{tikz}
\usetikzlibrary{decorations.pathreplacing,calligraphy}

\usepackage{algorithm}
\usepackage{algpseudocode}
\usepackage[square, comma, sort&compress,numbers]{natbib}
\usepackage{float}

\newtheoremstyle{plainNoItalics}{}{}{\normalfont}{}{\bfseries}{.}{ }{}
\theoremstyle{plain}

\theoremstyle{plainNoItalics}

\newtheorem{example}{Example}[section]
\theoremstyle{definition}
\theoremstyle{definition}
\theoremstyle{definition}
\theoremstyle{definition}\newtheorem{theorem}{Theorem}
\theoremstyle{definition}
\theoremstyle{definition}\newtheorem{rem}{Remark}
\numberwithin{equation}{section}
\numberwithin{theorem}{section}
\numberwithin{lemma}{section}
\numberwithin{figure}{section}
\numberwithin{example}{section}
\numberwithin{rem}{section}

\usepackage{lineno}
\usepackage{xcolor}
\usepackage{natbib}

\newcommand{\bU}{{\bm{U}}}
\newcommand{\bv}{{\bm{v}}}
\newcommand{\bxi}{\bm{\xi}}
\newcommand{\bd}{{\bm{d}}}
\newcommand{\bW}{{\bm{W}}}

\graphicspath{{figures/}}

\title{A learned conservative semi-Lagrangian finite volume scheme for transport simulations}
\author{
Yongsheng Chen
\thanks{School of Mathematical Sciences, Zhejiang University, Hangzhou, 310027, China. {\tt 22035024@zju.edu.cn}}
\and
Wei Guo
\thanks{Corresponding author. Department of Mathematics and Statistics, Texas Tech University, Lubbock, TX, 70409, USA. 
{\tt weimath.guo@ttu.edu}. }
\and
Xinghui Zhong
\thanks{ School of Mathematical Sciences, Zhejiang University, Hangzhou, 310027, China. {\tt zhongxh@zju.edu.cn}}
}

\begin{document}

\maketitle
\begin{abstract}
Semi-Lagrangian (SL) schemes are known as a major numerical tool for solving transport equations with many advantages and have been widely deployed in the fields of computational fluid dynamics, plasma physics modeling, numerical weather prediction,  among others. In this work, we develop a novel machine learning-assisted approach to accelerate the conventional SL finite volume (FV) schemes. The proposed scheme avoids the expensive tracking of upstream cells but attempts to learn the SL discretization from the data by incorporating specific inductive biases in the neural network, significantly simplifying the algorithm implementation and leading to improved efficiency. In addition, the method delivers sharp shock transitions and a level of accuracy that would typically require a much finer grid with traditional transport solvers. Numerical tests demonstrate the effectiveness and efficiency of the proposed method.

\end{abstract}
\section{Introduction}

Transport phenomena are ubiquitous in nature and are characterized by a type of partial differential equations (PDEs), namely  transport equations. Efficiently simulating transport equations  represents a great challenge mainly due to the potential presence of multiple spatiotemporal scales and discontinuities in the solution structure. 
The last several decades have witnessed tremendous developments of numerical methods for transport equations, yielding numerous successful applications in both science and engineering. Among them, the semi-Lagrangian (SL) transport schemes  attract lots of attention \cite{courant1952solution}, and have become a
major tool in simulating transport equations arising from numerical weather prediction \cite{robert1981stable,staniforth1991semi} and plasma simulations \cite{cheng1976integration,sonnendrucker1999semi}. Such an approach updates the mesh-based solution of a transport equation  by tracing the characteristics, offering several computational benefits, including high order accuracy and unconditional stability. Additionally, by incorporating an advanced discretization technique such as the weighted essentially non-oscillatory (WENO) method \cite{jiang1996efficient} or the discontinuous Galerkin method \cite{cockburn2001runge} into the SL framework, the resulting methods, see e.g., \cite{qiu2010conservative,qiu2011conservative,rossmanith2011positivity,qiu2011positivity},  attain outstanding performance for transport simulations. On the other hand, despite the effectiveness of the SL methodology, the direct numerical simulations may require an exceedingly large amount of computational resources when full resolution of fine-scale structures is of interest. Another difficulty associated with the SL approach is that tracking the characteristics accurately, including tracing and approximating upstream cells, is not only expensive in computational cost but also demands complex algorithm implementation, especially in high dimensions. For instance, the SL method proposed in \cite{lauritzen2010conservative} implements a rather sophisticated search algorithm.

With the rapid advancements in computing power and machine learning (ML) software over the past few decades, integrating ML techniques in simulating PDEs has become a thriving area of research. 
One notable example  is  the physics informed neural networks (PINNs) \cite{raissi_physics_2017,raissi_physics_2017-1}, where the solution to the underlying PDE is parameterized with a neural network (NN) and trained using a physics informed loss function. The partial derivatives of the approximate solution are obtained via automatic differentiation \cite{baydin_2018_automatic}. As an effective alternative to the traditional numerical algorithms, PINNs have been widely deployed to  solve complex problems in various fields, see \cite{jagtap_conservative_2020,yu_gradient-enhanced_2022,lu_physics-informed_2021,mcclenny_self-adaptive_2022,raissi_physics-informed_2019,lu_deepxde_2021,pang_fpinns_2019}. A comprehensive review of the literature on PINNs was provided in  \cite{cuomo_scientific_2022}.  
We also mention one popular research direction of leveraging  NNs to improve overall performance for the  traditional numerical methods, and related works include the troubled-cell indicator in \cite{ray_artificial_2018},  weights estimation to enhance WENO scheme in \cite{wang_learning_2020}, shock detection for WENO schemes in \cite{sun_convolution_2020}, and the total variation bounded constant estimation for limiters in \cite{yu_multi-layer_2022}.

Another group of NN-based methods for simulating PDEs is the so-called neural operator approach.  Two widely recognized works are  the DeepONet \cite{lu_deeponet_2021} which consists of two sub-networks and learns the nonlinear operators for identifying differential
equations, and the Fourier neural operator (FNO) \cite{li_fourier_2021}
which is designed based on parameterization of the integral kernel in the Fourier space and most closely resembles the reduced basis method.
A performance comparison of these two  methods is documented in \cite{lu_comprehensive_2022}. Other works in this approach include but not limited to \cite{kissas_learning_2022,li_neural_2020,bhattacharya_model_2021,trask_gmls-nets_2019}. 
%
For time-dependent PDEs, the cell-averaged neural network (CANN) was developed in \cite{qiu_cell-average_2021,chen_cell-average_2022}, which explores the approximation of the cell average difference of the solution between two consecutive time steps. Another different approach known as autoregressive methods was developed in \cite{bsl2019,greenfeld2019learning,hsieh2019learning,brandstetter2022message}, where the PDEs are  simulated iteratively, resembling conventional numerical methods with time marching. Specifically, an autoregressive method predicts the solution at the next time step with a NN based on the current state, and this approach is closely related to flow map learning methods; see, e.g., \cite{chen2022deep,churchill_deep_2022}.


In this paper, we consider the conservative SL finite volume (FV) formulation proposed in \cite{lauritzen2010conservative}, and  propose a novel ML-based approach to enhance the performance and accelerate the computation. The proposed methodology belongs to the category of autoregressive methods and  is motivated by the recent ML-based discretization for PDEs \cite{bar2019learning,zhuang2021learned,kochkov_machine_2021}. In contrast to the neural operator methods, the methods proposed in \cite{bar2019learning,zhuang2021learned,kochkov_machine_2021} are designed under the classical finite difference or FV framework, while a key component is replaced by an ML technique for improved performance. In particular,  instead of  polynomials, such an approach employs NNs to learn an optimal discretization for approximating derivatives. Additionally, the architecture of such ML-based methods can benefit from Tensor Processing Units to accelerate computation \cite{kochkov_machine_2021}, hence substantially reducing the computational cost.

Instead of the Eulerian method-of-line framework employed in \cite{bar2019learning,zhuang2021learned,kochkov_machine_2021}, our method is built in the SL setting, and the discretization is learned using a convolutional NN (CNN) architecture \cite{lecun1995convolutional,lecun1998gradient} and incorporating specific  inductive biases. In particular, we explicitly include a key quantity called the normalized shifts as part of the input in the NN, observing that such quantities contributes to computing the traditional SL FV discretization but in a rather complex manner, see Section \ref{sec:algorithm} below. With high-resolution high-fidelity data, the proposed method can effectively learn and predict the SL discretization. Hence, by replacing the most expensive and complex component of the SL formulation with an ML-based approach, the proposed method avoids the explicit implementation in tracing upstream cells as in \cite{lauritzen2010conservative}, achieving enhanced efficiency.  In addition,   as with the ML-based method using CNNs \cite{zhuang_learned_2021}, the learned SL discretization can accurately capture fine-scale features of the solution of interest with a coarse grid which often demands much higher grid resolution for a traditional discretization, leading to significant computational savings. We further propose to add a constraint layer to the NN to enforce exact mass conservation, which is a highly desired property for transport simulation and plays a vital role in long term accuracy and stability. By incorporating such an additional level of inductive bias, the ML model can enjoy improved generalization.  

It follows from the dependence of  the proposed methodology  on the CNNs that the discretizations are restricted to a Cartesian grid. 
Different from the traditional SL formulations, the proposed ML-based method in this paper employs a set of fixed stencils for the SL reconstruction, which incurs the time step restriction due to the Courant-Friedrichs-Lewy (CFL) condition for numerical stability.  However, it is observed in 
the numerical experiments in Section \ref{sec:num} that the CFL number  can be set as large as 1.8 if a fixed stencil consisting of five cells is employed, which is often adequate and efficient for many applications. It is worth mentioning that, although
the traditional SL approach is unconditionally stable in theory,
choosing an exceedingly large CFL increases the complexity of the search algorithm for tracing and approximating upstream cells, especially in high dimensions, and also hinders parallel efficiency; see \cite{lauritzen2010conservative,guo2014conservative}.


The rest of the paper is organized as follows. In Section \ref{sec:algorithm}, we first review the conventional SL FV scheme for 1D and two-dimensional (2D) transport equations, and provide a sufficient condition for mass conservation. Then we introduce the proposed conservative ML-based SL FV method based on a CNN architecture. In Section \ref{sec:num}, numerical results are presented to demonstrate the performance of the proposed method. The conclusion and future works are discussed in Section \ref{sec:con}.

\section{Algorithm}
\label{sec:algorithm}
\subsection{Semi-Lagrangian finite volume scheme}
\begin{figure}[h!]
	\centering
	\begin{tikzpicture}[x=1cm,y=1cm]
		\begin{scope}[thick]
			\draw[black]  (-3,0) -- (3,0)
			node[right]{$t^{m}$};
			\draw[black] (-3,2) -- (3,2)
			node[right]{$t^{m+1}$};
			
			
			\draw[black] (-2.4,0) -- (0,0);
			\draw[black] (0,0) -- (2.4,0);
			\draw[black] (0,2) -- (0,2.2)node[above] {$x_{i-\frac12}$};
			\draw[black] (0,0) -- (0,0.2);
			\draw[black] (2.4,2) -- (2.4,2.2)node[above] {$x_{i+\frac12}$};
			\draw[black] (2.4,0) -- (2.4,0.2);
			\draw[black] (-2.4,2) -- (-2.4,2.2);
			\draw[black] (-2.4,0) -- (-2.4,0.2);            
			\draw[-latex,dashed](0,2) to[out=200,in=70] (-1.4,0) node[above left=1pt] {$\tilde{x}_{i-\frac12}$};  
			\draw[-latex,dashed](2.4,2) to[out=200,in=70] ( 1.,0) node[above left=1pt] {$\tilde{x}_{i+\frac12}$};
			\draw [decorate,decoration={brace,amplitude=5pt,mirror,raise=0.5ex}]
			(-1.4,0) -- (1,0)node[midway,yshift=-1.5em] {$\tilde{I}_{i}$};
			
		\end{scope}
	\end{tikzpicture}

	\caption{Schematic illustration of the 1D SL FV scheme.\label{fig:slfv}}	
\end{figure}
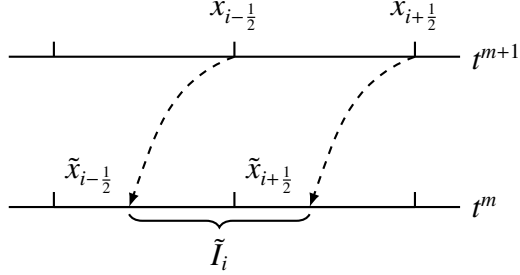

In this section, we review the classical SL FV scheme with remapping for linear transport equations; see e.g., \cite{lauritzen2010conservative,erath2013mass}. We start with the following one-dimensional (1D) equation in the conservative form 
\begin{equation}
	\label{eq:transport1d}
	u_t + (a(x,t)u)_x = 0,\quad x\in\Omega,
\end{equation}
where $a(x,t)$ is the velocity function. For simplicity, consider a uniform partition of the domain $\Omega$ with $N$ cells, i.e., $\Omega=\bigcup_{i=1}^N I_i$, where $I_i=[x_{i-\frac12},x_{i+\frac12}]$. Denote the mesh size as $h=|I_i|$.  The numerical solution  $\{U_i^m\}_{i=1}^N$ approximates the cell averages at the time step $t^m$. The SL FV method updates the cell averages to next time step $t^{m+1}$ by tracking the characteristic lines governed by the ordinary differential equation
\begin{equation}
	\label{eq:characteristic}
	\frac{dx(t)}{dt} = a(x(t),t).
\end{equation}
The upstream cell $\tilde{I}_i=[\tilde{x}_{i-\frac12},\tilde{x}_{i+\frac12}]$ of $I_i=[x_{i-\frac12},x_{i+\frac12}]$ is defined by evolving \eqref{eq:characteristic} with final values $x(t^{m+1}) = x_{i\pm\frac12}$ backward to $t^m$, as shown in Figure \ref{fig:slfv}. 

Then the exact solution $u(x,t)$ satisfies 
\begin{equation}
	\label{eq:slexact}
	\int_{I_i} u(x,t^{m+1}) dx = \int_{\tilde{I}_i} u(x,t^m) dx,\quad i=1,\ldots,N,
\end{equation}
based on which we can formulate the SL FV method with remapping as follows.  In each cell $I_j$, we reconstruct a polynomial $\phi_{j}(x)$ of degree $K$ using the $K+1$ cell averages from its neighboring cells in the spirit of the standard FV methodology, aiming for $(K+1)$-th order accuracy. Then, the right hand side of \eqref{eq:slexact} can be approximated by integrating the reconstructed polynomials in  $\tilde{I}_i$. Meanwhile, it is noted that reconstructed polynomials are defined piece-wisely  and hence are discontinuous across the cell interfaces. In addition, $\tilde{I}_i$ will not coincide with a background Eulerian cell in general. Therefore, the integrals must be computed in a subcell-by-subcell fashion. In particular, we denote by  $\Tilde{I}_{i,(j)}$ as the intersection of  $\tilde{I}_i$ and the background Eulerian cell $I_j$. The collection of indices of the Eulerian cells which have non-empty intersection with the upstream cell $\tilde{I}_i$ is defined as
$$
\mathcal{L}_i = \{j:\tilde{I}_i\cap I_j \ne \emptyset\}.
$$
We write $\phi_j(x)=\sum_{k=0}^K c_j^{(k)}v_j^{(k)}(x)$, where $v_j^{(k)}(x)$ are the local polynomial basis in the cell $I_j$. At the discrete level, \eqref{eq:slexact} becomes 
\begin{equation}
	\label{eq:sl}
	U_i^{m+1}|I_i| = \sum_{j\in\mathcal{L}_i}\sum_{k=0}^K c_{j}^{(k)} \int_{\tilde{I}_{i,(j)}} v_j^{(k)}(x)dx \doteq  \sum_{j\in\mathcal{L}_i}\sum_{k=0}^K c_{j}^{(k)}\omega_{i,j}^{(k)},
\end{equation}
yielding the SL FV scheme of order $K+1$.  A key observation of the SL FV scheme is that $\{U_i^{m+1}\}$ are fully determined by cell averages at previous time step $\{U_i^{m}\}$ together with the normalized shifts 
$$\xi_{i-\frac12} := \frac{\tilde{x}_{i-\frac12}-{x}_{i-\frac12}}{h},\quad i=1,\ldots,N,$$
which encode the geometry information of the upstream cells. In particular, the coefficients $\{\omega_{i,j}^{(k)}\}$ in \eqref{eq:sl} can be expressed using  $\{\xi_{i-\frac12}\}$. More detains can be found in \cite{qiu2011conservative}.
Therefore, the SL FV scheme \eqref{eq:sl} can be recast into the following formulation
\begin{equation}
	\label{eq:sl_re}
	U_i^{m+1} = \sum_{\ell\in\mathcal{S}^m_i} d^m_{i,\ell} U_{\ell}^m,
\end{equation}
where $\mathcal{S}^m_i$ denotes the stencil employed to update  $U_i^{m+1}$, and $\{d_{i,\ell}^m\}$ are the associated coefficients. It is worth mentioning that $\{d_{i,\ell}^m\}$ may depend on the numerical solution if a nonlinear reconstruction is used.


It can be shown that the SL FV scheme \eqref{eq:sl} conserves mass in the discrete sence, i.e.,
\begin{equation}
	\label{eq:conservative}
	\sum_{i} U_i^{m+1} = \sum_{i} U_i^{m}.
\end{equation}
For an FV method in the form of \eqref{eq:sl_re}, the  conditions for mass conservation is discussed in the following theorem.
\begin{theorem}\label{thm:conservation} An FV method in the form of \eqref{eq:sl_re} is mass conservative, i.e., satisfying \eqref{eq:conservative}, when
	\begin{equation}\label{eq:conser_condition}
		\sum_{i \in \mathcal{E}_\ell^m} d_{i,\ell}^m   =1,\quad\forall l,
	\end{equation}
	where  $\mathcal{E}_\ell^m =\{i: \ell\in \mathcal{S}_i^m\}
	$
	is the index set for a collection of cells $\{I_i\}$ for which $U_\ell^m$  contributes to the update of $\{U_i^{m+1}\}$, i.e., the region of influence of $U_\ell^m$.
	
	If the FV method is linear, i.e., the coefficients $d^m_{i,\ell}$ is independent of the numerical solution $\{U_i^m\}$, then \eqref{eq:conser_condition} is also a necessary condition of mass conservation.
	
\end{theorem}

\noindent{\bf Proof.}  Summing \eqref{eq:sl_re} over $i$ gives
\begin{align*}
	\sum_i U_i^{m+1} &= \sum_i\sum_{\ell\in\mathcal{S}^m_i} d^m_{i,\ell} U_{\ell}^m\\
	&=\sum_{\ell}\left(\sum_{i \in \mathcal{E}_\ell^m} d_{i,\ell}^m\right) U_\ell^m.
\end{align*}
Then, if \eqref{eq:conser_condition} holds, then $\sum_i U_i^{m+1}=\sum_\ell U_\ell^{m}$, and hence the scheme is conservative. 

Further, if the scheme is linear and conservative, then 
$$
\sum_{\ell}\left(\sum_{i \in \mathcal{E}_\ell^m} d_{i,\ell}^m\right) U_\ell^m = \sum_{\ell} U_\ell^m
$$
holds for arbitrary $\{U_\ell^m\}$. By letting $U_\ell^m = 1$ and $U_j^m=0$ for $j\ne\ell$, we have \eqref{eq:conser_condition}. \qedsymbol

Hence, \eqref{eq:conser_condition}  provides a feasible  way to enforce mass conservation for any FV schemes expressed in the form of  \eqref{eq:sl_re}. 

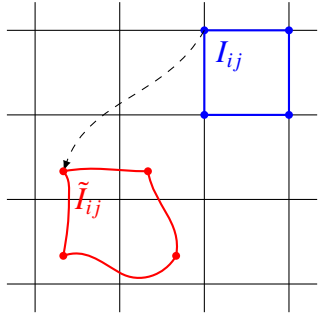
\begin{figure}[h]
	\centering
	\begin{tikzpicture}[scale = 0.75]
		\draw[black,thin] (0,0.5) node[left] {} -- (5.5,0.5)
		node[right]{};
		\draw[black,thin] (0,2.) node[left] {$$} -- (5.5,2)
		node[right]{};
		\draw[black,thin] (0,3.5) node[left] {$$} -- (5.5,3.5)
		node[right]{};
		\draw[black,thin] (0,5 ) node[left] {$$} -- (5.5,5)
		node[right]{};
		\draw[black,thin] (0.5,0) node[left] {} -- (0.5,5.5)
		node[right]{};
		\draw[black,thin] (2,0) node[left] {$$} -- (2,5.5)
		node[right]{};
		\draw[black,thin] (3.5,0) node[left] {$$} -- (3.5,5.5)
		node[right]{};
		\draw[black,thin] (5,0) node[left] {$$} -- (5,5.5)
		node[right]{};
		\fill [blue] (3.5,3.5) circle (2pt) node[] {};
		\fill [blue] (5,3.5) circle (2pt) node[] {};
		\fill [blue] (3.5,5) circle (2pt) node[below right] {$I_{ij}$} node[above left] {};
		\fill [blue] (5,5) circle (2pt) node[] {};
		
		\draw[thick,blue] (3.5,3.5) node[left] {} -- (3.5,5)
		node[right]{};
		\draw[thick,blue] (3.5,3.5) node[left] {} -- (5,3.5)
		node[right]{};
		\draw[thick,blue] (3.5,5) node[left] {} -- (5,5)
		node[right]{};
		\draw[thick,blue] (5,3.5) node[left] {} -- (5,5)
		node[right]{};
		\draw[-latex,dashed](3.5,5)node[right,scale=1.0]{}
		to[out=240,in=70] (1,2.50) node[] {};
		
		\draw (0.5+0.01,2-0.01) node[fill=white,below right] {};
		
		\draw [red,thick] (1,1)node[right,scale=1.0]{}
		to[out=20,in=150] (2,0.7) node[] {};
		
		\draw [red,thick] (2,0.7)node[right,scale=1.0]{}
		to[out=330,in=240] (3,1) node[] {};
		\draw [red,thick] (1,2.5)node[right,scale=1.0]{}
		to[out=310,in=90] (1.1,2) node[] {};
		\draw [red,thick] (1.1,2)node[right,scale=1.0]{}
		to[out=270,in=80] (1,1) node[] {};

		\draw [red,thick] (1,2.5)node[right,scale=1.0]{}
		to[out=10,in=180] (2.5,2.5) node[] {};
		
		\draw [red,thick] (3,1)node[right,scale=1.0]{}
		to[out=80,in=280] (2.5,2.5) node[] {};
		\fill [red] (1.,1) circle (2pt) node[above right,black] {};
		\fill [red] (3,1) circle (2pt) node[] {};
		\fill [red] (1,2.5) circle (2pt) node[below right] {$\tilde{I}_{ij}$} node[above left] {};
		\fill [red] (2.5,2.5) circle (2pt) node[] {};
	\end{tikzpicture}
	\caption{Schematic illustration of the 2D SL FV scheme.\label{fig:slfv2d}}	
\end{figure}

The FV SL scheme with remapping can be extended to the 2D linear transport equation with variable coefficients
\begin{equation}\label{eq:2d}
	u_t + \nabla\cdot (\bv(x,y,t) u) = 0,\quad (x,y)\in\Omega, 
\end{equation}
where $\bv$ denotes the velocity field $\bv(x,y,t) = (a(x,y,t),b(x,y,t))$. The associated characteristic system writes
\begin{equation}
	\label{eq:characteristic2d}
	\begin{cases}
		\frac{dx(t)}{dt} &= a(x(t),y(t),t),\\
		\frac{dy(t)}{dt} &= b(x(t),y(t),t).\\
	\end{cases}
\end{equation}
Assume the domain $\Omega$ is partitioned uniformly with a collection of rectangular cells, i.e., $\Omega=\bigcup_{i,j} I_{ij}$, where $I_{ij}=[x_{i-\frac12},x_{i+\frac12}]\times [y_{j-\frac12},y_{j+\frac12}]$. Denote by $h_x=x_{i+\frac12}-x_{i-\frac12}$ and $h_y=y_{j+\frac12}-y_{j-\frac12}$ the mesh sizes in $x$ and $y$ directions, respectively.
Similar to the 1D case, by solving \eqref{eq:characteristic2d} backward in time from $t^{m+1}$ to $t^m$, we obtain the upstream cell $\Tilde{I}_{ij}$ of cell $I_{ij}$, as shown in Figure \ref{fig:slfv2d}. Denote the upstream point of each grid point $(x_{i-\frac12},y_{j-\frac12})$ as $(\tilde{x}_{i-\frac12,j-\frac12},\tilde{y}_{i-\frac12,j-\frac12})$, and then
define 
$$
\xi_{i-\frac12,j-\frac12} =\frac{\tilde{x}_{i-\frac12,j-\frac12} - x_{i-\frac12}}{h_x},\quad \eta_{i-\frac12,j-\frac12} =\frac{\tilde{y}_{i-\frac12,j-\frac12} - y_{j-\frac12}}{h_y}
$$
as the normalized shifts of grid point $(x_{i-\frac12},y_{j-\frac12})$ in the $x$ and $y$ directions, respectively.

The 2D SL FV scheme is then formulated with the following identity \begin{equation}
	\label{eq:slexact2d}
	\iint_{I_{ij}} u(x,y,t^{m+1}) dxdy = \iint_{\tilde{I}_{ij}} u(x,y,t^m) dxdy.
\end{equation} 
Denote by $U^m_{ij}$ the cell average of the numerical solution in the cell $I_{ij}$ at time step $t=t^m$. Similar to the 1D case, to update the numerical solution,  a polynomial is reconstructed over each cell using cell averages, and the integral on the right-hand side of \eqref{eq:slexact2d} is computed in the subcell-by-subcell fashion.  We write the scheme as
\begin{equation}
	\label{eq:sl_re2d}
	U_{ij}^{m+1} = \sum_{\ell\in\mathcal{S}^m_{ij}} d^m_{ij,\ell} U_{\ell}^m,
\end{equation}
where $\mathcal{S}^m_{ij}$ denotes the stencil employed to update  $U_{ij}^{m+1}$. Again, the coefficients $\{d_{ij,\ell}^m\}$ are determined by the solution averages $\{U_{ij}^m\}$ together with normalized shifts $\{\xi_{i-\frac12,j-\frac12}^m\}$, $\{\eta_{i-\frac12,j-\frac12}^m\}$. In addition, as with the 1D case, it can be shown that the scheme \eqref{eq:sl_re2d} is mass conservative if 
\begin{equation}\label{eq:conser_condition2d}
	\sum_{(i,j) \in \mathcal{E}_\ell^m} d_{ij,\ell}^m   =1,\quad\forall l,
\end{equation}
where $\mathcal{E}_\ell^m =\{(i,j): \ell\in \mathcal{S}_{ij}^m\}.$

\begin{rem}
	The SL FV methods are high order accurate in space but exact in time if the upstream cells are traced exactly. Further, the schemes are free of the CFL time step restriction for stability as the reconstructions are local.
	
\end{rem}

\begin{rem}
	The most computationally intensive part of the SL FV methods introduced above lies in tracking the geometry information of the upstream cells, including organizing the overlap regimes and integrating the local polynomial basis in a subregime-by-subregime manner. Such a search algorithm becomes highly challenging and expensive in high dimensions when the upstream cells may deform into irregular shapes as shown in Figure \ref{fig:slfv2d}. In next section, we propose an ML-assisted SL FV method to avoid such expensive  tracing of upstream cells.
\end{rem}

\subsection{Data-driven conservative SL FV scheme}
In this section, we introduce a novel data-driven SL FV scheme with enhanced accuracy and efficiency. This is motivated by a class of successful ML-based approaches for optimal discretizations for PDEs \cite{bar2019learning,zhuang2021learned,kochkov_machine_2021}. The main idea of such methods is that the solution manifold of a PDE often exhibits low dimensional structures, such as the recurrent patterns and  coherent structures. Given high resolution training  data, the ML-based discretization can effectively parameterize the solution manifold with coarse grids and attain a level of accuracy which often requires an order-of-magnitude finer grid for a  traditional method using polynomial-based approximations. Hence, the methods can potentially capture the dynamics of interest even with an under-resolved grid.  Furthermore, under the standard method-of-line framework, it is natural for the methods to satisfy inherent physical constraints such as conservation of mass, momentum and energy of the underlying physical system, as opposed to the purely data driven approach, and such an inductive bias aids in reliability and generalization of the ML-based model. The proposed data-driven SL FV scheme aims to take advantage of the methodology of ML-based discretization \cite{bar2019learning,zhuang2021learned} and the SL FV formulation reviewed in the previous section for efficient transport simulations with the mass conservation.

\begin{figure}[!htbp]
	\centerline{\includegraphics[width=0.98\textwidth]{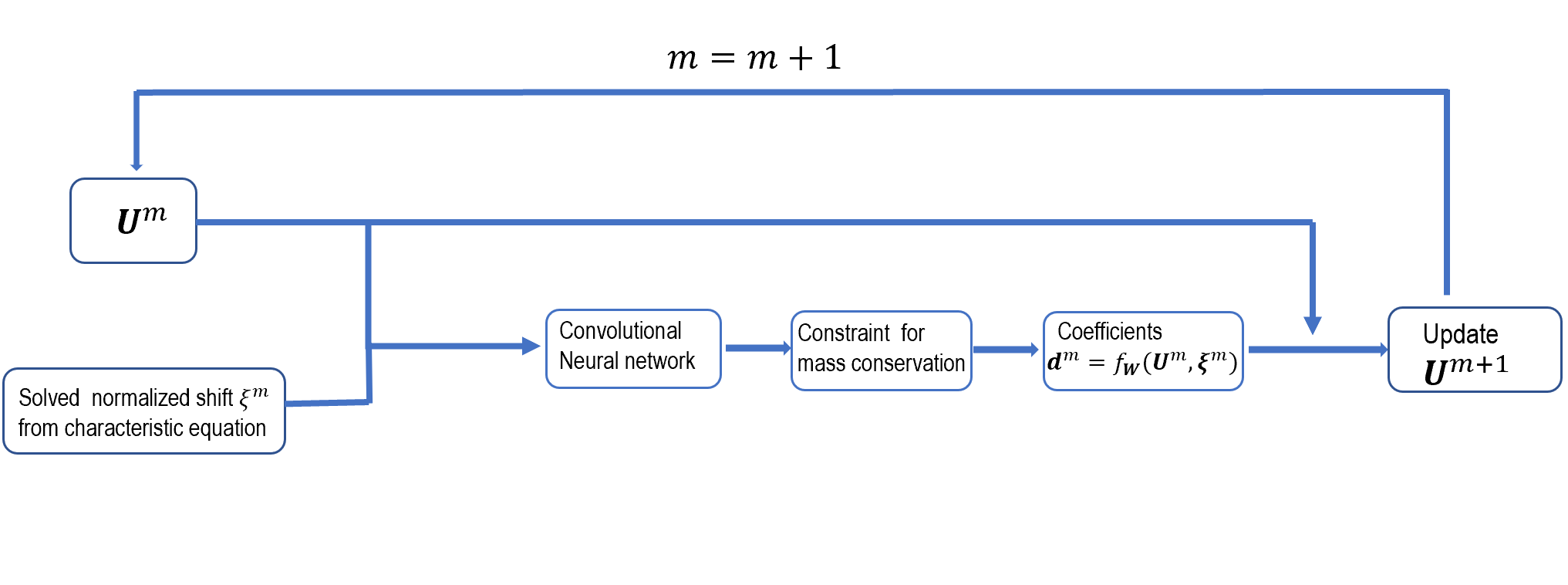}}
	\caption{Illustration of the proposed ML-based SL FV method.}\label{fig:network_structure}
\end{figure}

We  first consider the 1D case and then briely discussed the generalization to the 2D case. Without abuse of notations, we denote by $\bU^m$, $\bxi^m$, and $\bd^m$ the collection of $\{U^m_i\}$, $\{\xi^m_{i-\frac12}\}$, and $\{d^m_{i,\ell}\}$, respectively. The proposed ML-based SL FV method is schematically illustrated in Figure \ref{fig:network_structure}.

The main idea is that instead of employing formulation \eqref{eq:sl}, we 
propose to work with \eqref{eq:sl_re} in which the coefficients are inferred by a trained feed-forward NN.    Observing that $\bd^m$ is determined by $\bU^m$ and $\bxi^m$ for the SL FV scheme, we design 
\begin{equation}
	\label{eq:fw}
	\bd^m = f_\bW(\bU^m,\bxi^m),
\end{equation} 
where the NN $f_\bW$ takes $\bU^m$ and  $\bxi^m$ as two-channel input and is constructed as a stack of convolutional layers with trainable parameters $\bW$ and nonlinear activation functions, such as ReLU. The proposed method attempts to replace the most expensive component of the traditional SL methods with a data-driven approach, significantly simplifying the algorithm implementation and improving efficiency. In addition, the employed CNN can effectively extract hierarchical features of the solution and meanwhile enable translation equivalence \cite{lecun1989handwritten,lecun1995convolutional}, which are highly desired for transport modelling. With high-resolution training data, the proposed model is capable of learning the optimal coefficients for SL transport. It is numerically observed in Section \ref{sec:num} that the proposed ML-based SL FV method  outperforms immensely the famous WENO method \cite{jiang1996efficient} for a collection of benchmark tests.

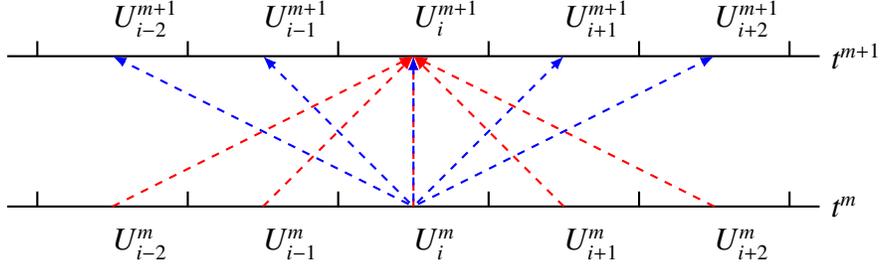
\begin{figure}[h!]
	\centering
	\begin{tikzpicture}[x=1cm,y=1cm]
		\begin{scope}[thick]
			\draw[black]  (-5.4,0) -- (5.4,0)
			node[right]{$t^{m}$};
			\draw[black] (-5.4,2) -- (5.4,2)
			node[right]{$t^{m+1}$};
			

			\draw[black] (-5,2) -- (-5,2.2);
			\draw[black] (-3,2) -- (-3,2.2);
			\draw[black] (-1,2) -- (-1,2.2);
			\draw[black] (1,2) -- (1,2.2);
			\draw[black] (3,2) -- (3,2.2);
			\draw[black] (5,2) -- (5,2.2);
			\draw[black] (-5,0) -- (-5,0.2);
			\draw[black] (-3,0) -- (-3,0.2);
			\draw[black] (-1,0) -- (-1,0.2);
			\draw[black] (1,0) -- (1,0.2);
			\draw[black] (3,0) -- (3,0.2);
			\draw[black] (5,0) -- (5,0.2);
			\node[text width=0cm] at (0,-0.5) 
			{$U^m_{i}$};
			\node[text width=0cm] at (2,-0.5) 
			{$U^m_{i+1}$};
			\node[text width=0cm] at (4,-0.5) 
			{$U^m_{i+2}$};
			\node[text width=0cm] at (-2,-0.5) 
			{$U^m_{i-1}$};
			\node[text width=0cm] at (-4,-0.5) 
			{$U^m_{i-2}$};
			\node[text width=0cm] at (0,2.5) 
			{$U^{m+1}_{i}$};
			\node[text width=0cm] at (2,2.5) 
			{$U^{m+1}_{i+1}$};
			\node[text width=0cm] at (4,2.5) 
			{$U^{m+1}_{i+2}$};
			\node[text width=0cm] at (-2,2.5) 
			{$U^{m+1}_{i-1}$};
			\node[text width=0cm] at (-4,2.5) 
			{$U^{m+1}_{i-2}$};
			\draw[red,-latex,dashed]  (-4,0) to (0,2) ;  
			\draw[red,-latex,dashed] (-2,0) to (0,2) ;
			\draw[red,-latex,dashed](0,0) to (0,2) ;
			\draw[red,-latex,dashed](2,0) to (0,2) ;
			\draw[red,-latex,dashed] (4,0) to (0,2) ;
			
			\draw[blue,-latex,dashed,dash phase=4pt]  (0,0) to (-4,2) ;  
			\draw[blue,-latex,dashed,dash phase=4pt] (0,0) to (-2,2) ;
			\draw[blue,-latex,dashed,dash phase=4pt](0,0) to (0,2) ;
			\draw[blue,-latex,dashed,dash phase=4pt](0,0) to (2,2) ;
			\draw[blue,-latex,dashed,dash phase=4pt] (0,0) to (4,2) ;
			
		\end{scope}
	\end{tikzpicture}

	\caption{Schematic illustration of the FV SL scheme with  fixed 5-cell stencil $\{I_{i-2},I_{i-1},I_{i},I_{i+1},I_{i+2}\}$. $U_i^{m+1}$ directly depends on $\{U^m_{i-2},U^m_{i-1},U^m_{i},U^m_{i+1},U^m_{i+2}\}$ and $U^m_{i}$ contributes to computing $\{U^{m+1}_{i-2},U^{m+1}_{i-1},U^{m+1}_{i},U^{m+1}_{i+1},U^{m+1}_{i+2}\}$.\label{fig:slfv_fixed}}	
\end{figure}

Once $\bm{d}^m$ is obtained, the solution $\bU^{m+1}$ is updated with \eqref{eq:sl_re}.  However, unlike the standard SL formulation, we employ a set of fixed centered stencils $\mathcal{S}^m_{i} = \{i-s,i-s+1,\ldots,i+s-1,i+s\}$ containing $2s+1$ cells. Not only will this greatly simplify the algorithm development but also make it convenient to satisfy the mass conservation constraint.  For example, when $s=2$, we have the 5-cell stencil. The NN $f_{\bW}$ predicts the coefficients $\{d^m_{i,\ell},\,\ell=i-2,\ldots,i+2\}$, and $U^{m+1}_i$ is given by $U^{m+1}_i=\sum_{\ell =i -2}^{i+2} d^m_{i,\ell} U_\ell^m, $
as shown in Figure \ref{fig:slfv_fixed} . Unfortunately, using fixed stencils will destroy the desired unconditional stability of the SL method.  In \cite{larios2022error}, an ML-based SL approach is developed in the context of the level-set method. Such a method is designed based on correcting the local error incurred by the standard SL method using a NN and  resembles a localized version of the error-correcting method using ML \cite{pathak2020using}, and hence the CFL time step restriction can be avoided.

Additionally, mass conservation is known as a critical requirement for transport simulations, as it is directly related to the long term accuracy and stability. Hence, it is highly desired to build such an inductive bias of mass conservation into the ML model for better generalization. Fortunately, in light of Theorem \ref{thm:conservation}, it is straightforward to enforce mass conservation for the proposed ML model: since $U_{i}^m$  contributes to computing  $\{U_{i-s}^{m+1},U_{i-s+1}^{m+1},\ldots,U_{i+s}^{m+1}\}$, as long as $\sum_{j=i-s}^{i+s}d^m_{j,i}=1$, the mass is conserved, see Figure \ref{fig:slfv_fixed}. Hence, to enforce the condition, we simply add a constraint layer to $f_\bW$ in \eqref{eq:fw} for exact mass conservation. Note that it is nontrivial to enforce mass conservation with a non-fixed stencil under the employed CNN framework. 

The proposed algorithm can be generalized to the 2D case with a few simple modifications. We define the following NN
\begin{equation}
	\label{eq:nn2d}
	\bm{d}^m = f_{\bm{W}}(\bU^m,\bm{\xi}^m,\bm{\eta}^m),
\end{equation}
where, as with the 1D case, $\bU^m,\bm{\xi}^m,\bm{\eta}^m$ are 2D tensors defined by collecting $\{U_{ij}^m\}$, $\{\xi_{i-\frac12,j-\frac12}^m\}$, $\{\eta_{i-\frac12,j-\frac12}^m\}$, respectively. $f_{\bW}$ takes a 3-channel input and outputs the coefficient tensor $\bm{d}^m$ which is used to update the solution. The NN $f_{\bW}$ is constructed by staking a sequence of 2D convolutional layers together a constraint layer to enforce the condition \eqref{eq:conser_condition2d} for exact mass conservation. Last, note that similar to the 1D case, we employ a set of squared fixed stencils, each of which consists $(2s+1)\times(2s+1)$  cells. Hence, the SL formulation is constrained by the CFL condition.

Note that compared to the ML-assisted Eulerian transport method in \cite{zhuang_learned_2021}, the proposed SL FV method has the advantage of avoiding the use of any explicit time integrators. As a result, only a single evaluation of the neural network $f_\bW$ is required per time evolution. In addition, the incorporation of characteristic information allows for a higher CFL of up to 2 while maintaining numerical stability with the use of fixed 5-cell stencils (see Figure \ref{fig:slfv_fixed}).    

\section{Numerical results}
\label{sec:num}
In this section, we carry out a series of numerical experiments to demonstrate the performance of our data-driven SL FV scheme for a collection of 1D and 2D transport equations.  Noteworthy,
the performance of the proposed scheme depends on the choice of hyperparameters of the CNN structure, and  numerical results with default settings 
are presented in this section for simplicity.
For both 1D and 2D equations, we employ 6 convolutional layers with 32 filters per layer, utilizing a kernel size of 3 or 5 and $3\times3$ or $5\times5$ for 1D and 2D cases, respectively, similar to \cite{zhuang_learned_2021}. A constraint layer is added to ensure mass conservation. The activation function used is ReLU. Following \cite{Ilya_fix_2017}, Adam optimization algorithm is applied to train the network. We employ the Eulerian fifth-order FV WENO (WENO5) method \cite{jiang1996efficient}, combined with the third order strong-stability-preserving Runge-Kutta (SSPRK3) time integrator \cite{gottlieb2001strong} over fine-resolution grids to produce ground-truth reference solution trajectories, and the training data are generated by coarsening the reference solutions by a certain factor. For all the test examples, we mainly report the results by the proposed ML-based SL FV method and the WENO5 method with the same mesh resolution, together with the reference solutions for comparison. It is worth emphasizing that we can use any accurate and reliable transport methods (e.g. WENO5 + SSPRK3) to generate training data. In all the plots reported below, ``Neural net" denotes the proposed method and ``WENO5" denotes the WENO5 method combined with SSPRK3. 

\subsection{One-dimensional transport equations}
In this subsection, we present numerical results for simulating 1D transport equations.

\begin{example}\label{square}
	In this example we consider the following advection equation with a constant coefficient
	\begin{equation}\label{eq:const1d}
		u_t + u_x = 0,\quad x\in [0,1],
	\end{equation}
	and periodic conditions are imposed. 
\end{example}

The training data is generated by coarsening 30 high-resolution solution trajectories over a 256-cell grid by a factor of 8. The initial condition for each trajectory is   a   square wave with height randomly sampled from $[0.1,1]$ and width from $[0.2,0.4]$.
In addition, each coarsened trajectory contains 256 sequential time steps, and the CFL number is within the range of $[0.3,1.8]$. The centered $5$-cell stencils are employed to update the solution. For testing, initial conditions are square functions randomly selected from the same width and height range. For comparison, the reference solution is produced using WENO5  on a high-resolution 256-cell grid and then down-sampled to a coarse grid of 32 cells, the procedure of which is the same as generating the training data. 

\begin{figure}[!htbp]
	\centerline{\includegraphics[width=0.95\textwidth]{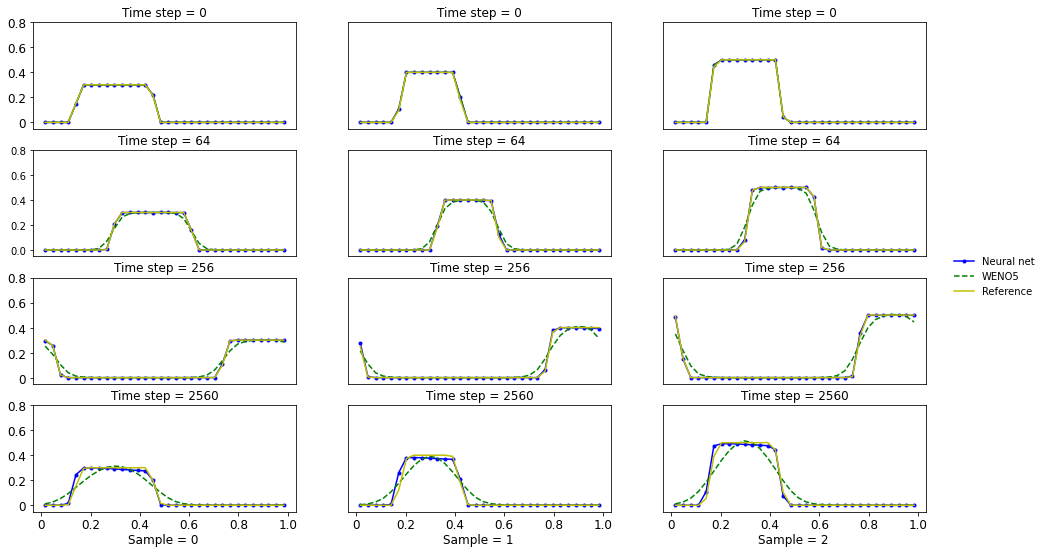}}
	\caption{Three test samples for square waves in Example \ref{square}.  CFL=0.6 for both our method and WENO5.}\label{fig:square_three_curve}
\end{figure}

Figure \ref{fig:square_three_curve} plots three test samples during forward integration at several instances of time with CFL = 0.6. It is observed that the WENO5 method exhibits significant smearing near discontinuities, which deteriorates over time as a result of the accumulation of numerical diffusion.
In contrast, the proposed ML-based solver has much improved shock resolution compared to WENO5  with the same mesh resolution: the numerical results are free of spurious oscillation and having very sharp shock transition. Notably, after long time simulations of 2560 time steps, which is 10 times of time steps for training, the results by the proposed method still stay very close to the reference solution, while a large amount of accumulated numerical diffusion of WENO5 dramatically smears the discontinuities. 

\begin{figure}[!htbp]
	\centering
	\subfigure[]{\includegraphics[width=0.45\textwidth]{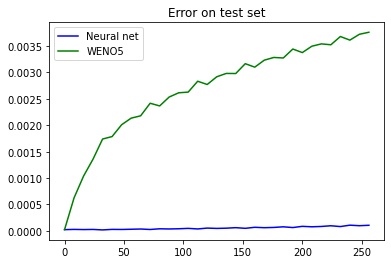}}
	\subfigure[]{\includegraphics[width=0.45\textwidth]{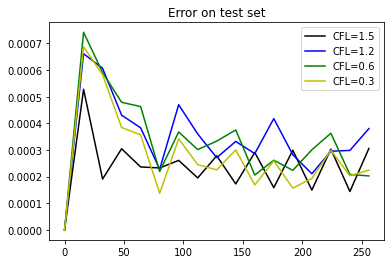}}
	\caption{Time histories of the errors in Example \ref{square}. (a): Errors by the proposed method  compared with WENO5 averaged over all test samples. (b): Errors by the proposed method with different CFLs.}\label{fig:square_error}
\end{figure}

\begin{figure}[!htbp]
	\centerline{\includegraphics[width=0.5\textwidth]{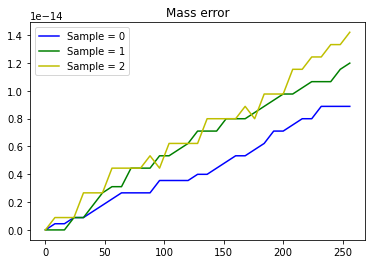}}
	\caption{Time histories of the  deviation of total mass for three test samples in Example \ref{square}. CFL=0.6.}\label{fig:square_mass}
\end{figure}

Figure \ref{fig:square_error} plots the time histories of the mean square errors which have been averaged over all test samples. Comparing with WENO5, see Figure \ref{fig:square_error}(a), the proposed ML-based method achieves a factor of approximately 8.6 less error in magnitude. Moreover, it is observed that the error by WENO5 increases over time, while the error by our method stays the same in magnitude with slight fluctuation. We further investigate the  performance of the proposed method  with different CFL numbers, i.e., different time step sizes, and report the result  in Figure \ref{fig:square_error}(b). We observe that the errors by our method with different CFLs are almost of the same magnitude over time.
Figure \ref{fig:square_mass} presents the time evolution of deviation in total mass for three test solution trajectories generated by our method. Evidently, the total mass is conserved up to the machine precision as expected.

\begin{figure}[!htbp]
	\centerline{\includegraphics[width=0.65\textwidth]{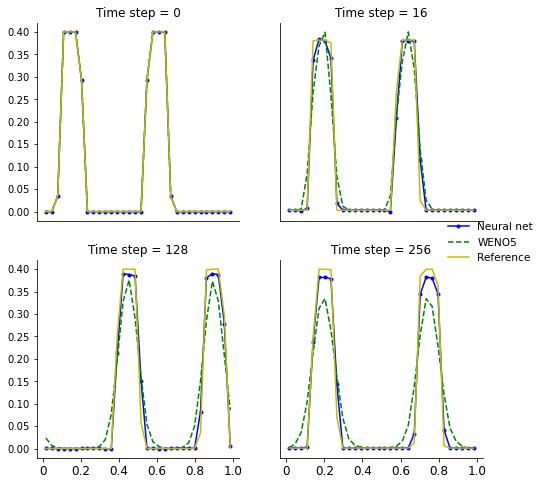}}
	\caption{Numerical solutions for Example \ref{square} with an initial condition of two square waves. The proposed ML model is trained with a data set for which each trajectory only contains a single square wave. CFL=0.6.}\label{fig:square_two}
\end{figure}

Although our solver is trained with a data set where each trajectory contains  a single square wave, it is observed that the model can be generalized to simulate  the advection equation \eqref{eq:const1d} with an initial condition consisting of two square waves, as shown in Figure \ref{fig:square_two}.  Again, our solver demonstrates superior performance over WENO5.  It is worth mentioning that  the traditional reduced order models employing a direct parameterization of the underlying solution manifold  are not capable of such generalization. 

Last, we remark that the CFL number cannot be chosen above 2, otherwise the loss would not decrease during training. Such an observation is partly attributed to the fact that the region of   dependence for updating one cell average is not completely contained within the 5-cell stencil if the CFL number is greater than 2.


\begin{example}\label{triangle}
	In this example, we consider the advection equation \eqref{eq:const1d} with a more complicated solution profile consisting of triangle and square waves. 
\end{example}

Similar to the previous example, we generate 30 high-resolution solution trajectories, each consisting of one triangle and one square waves with heights randomly sampled from $[0.2,0.8]$ and widths from $[0.2,0.3]$ over the 256-cell grid. By reducing the resolution from 256 to 32 (by a factor of 8), we obtain our training data. For this example, we set the stencil size to be $3$, and the maximum CFL number allowed for training is reduced to $0.975$ based on our numerical experiments.  Each solution trajectory in the training data set contains 256 time steps with the CFL number ranging in [0.3,0.975].  The test data are randomly sampled from the same ranges of width and height. To evaluate the performance of the proposed method, we calculate the ground-truth reference solution using WENO5 on a high-resolution grid with 256 cells and then reduce the resolution by a factor of 8 to a coarse grid of 32 cells. 

\begin{figure}[!htbp]
	\centerline{\includegraphics[width=0.95\textwidth]{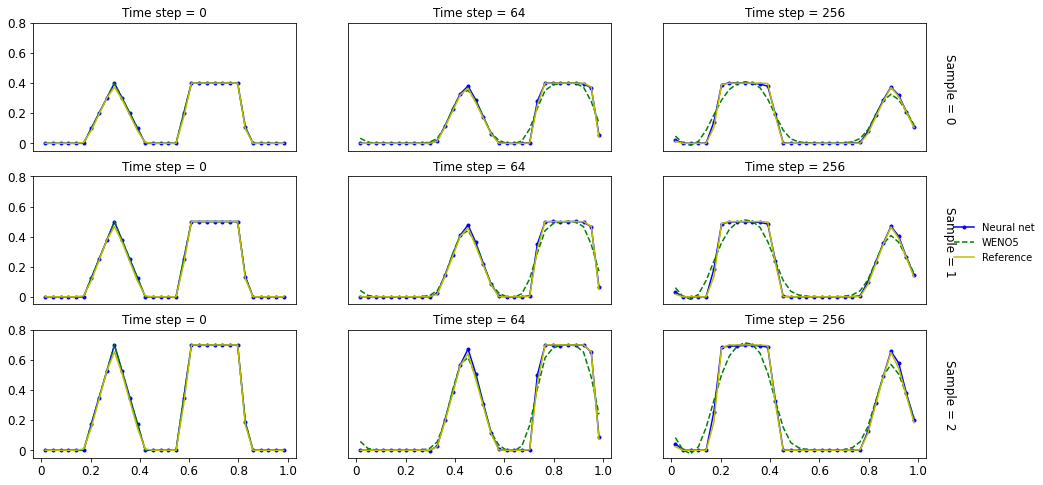}}
	\caption{Numerical solutions of three test samples for advection of triangle and square waves in Example \ref{triangle}.  CFL=0.6.}\label{fig:triangle_three_curve}
\end{figure}

\begin{figure}[!htbp]
	\centering
	\subfigure[]{\includegraphics[width=0.45\textwidth]{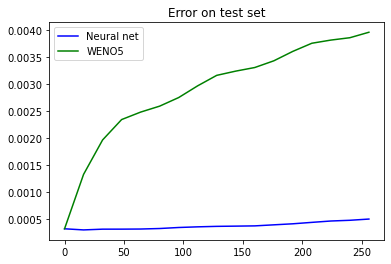}}
	\subfigure[]{\includegraphics[width=0.45\textwidth]{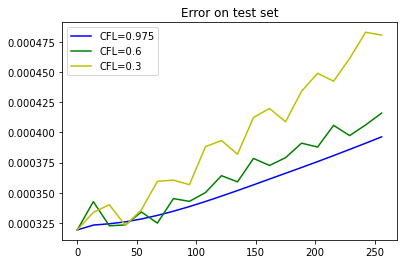}}
	\caption{Time histories of the errors in Example \ref{triangle}. (a): Errors by the proposed method  compared with WENO5 averaged over all test samples. (b): Errors by the proposed method with different CFLs.}
	\label{fig:triangle_error}
\end{figure}

\begin{figure}[!htbp]
	\centerline{\includegraphics[width=0.5\textwidth]{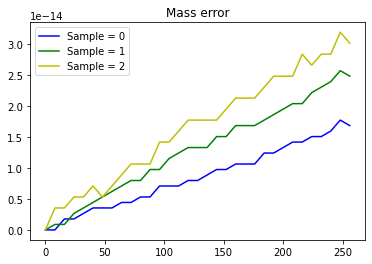}}
	\caption{Time evolution of the deviation of the total mass for three test samples in Exmpale \ref{triangle}. }\label{fig:triangle_mass}
\end{figure}

In Figure \ref{fig:triangle_three_curve}, 
we plot three test samples at several instances of time during forward integration with CFL = 0.6.  The proposed ML-based solver generates numerical results with significantly higher resolution of non-smooth structures compared to WENO5. We then present the time histories of the mean square error in Figure \ref{fig:triangle_error} (a), averaged over all test examples. Our ML-based solver achieves a reduction of error magnitude by a factor of approximately 7.8 compared to the traditional WENO5 solver. Additionally, the error of our solver grows at a much slower rate with time compared to that of WENO5. To further validate our solver, we consider three CFL numbers for testing, and present the time histories of errors in Figure \ref{fig:triangle_error}(b). It is observed that employing a larger CFL results in a smaller error and slower growth in time. Similar to the previous example, our method is mass conservative up to machine precision as demonstrated in Figure \ref{fig:triangle_mass}.

\begin{example}
	\label{eg:variable}
	We simulate the following 1D advection equation with a variable coefficient
	\begin{equation}\label{eq:variable-velocity}
		u_t+(u\sin(x+t))_x=0,\quad x\in[0,2\pi],
	\end{equation}
	and periodic conditions are imposed.
\end{example} 

We generate 90 solution trajectories, and each initial condition is a square function  with heights randomly sampled from $[0.1,1]$ and widths from $[2.5,3.5]$ over the high-resolution grid of 256 cells. By reducing the grid resolution with a factor of 8, we obtain the training data over a 32-cell grid. In addition, each trajectory in the training data contains 30 sequential time steps, with the CFL number ranging in   $[0.3,1.2]$. In addition, the center of each square function is randomly sampled from the whole domain $[0,2\pi]$. We set the stencil size to be $5$. Again, the reference solution is generated by WENO5 over the 256-cell grid  and down-sampled to the coarse grid of 32 cells. Note that the solution structure is more complicated than previous examples.

\begin{figure}[!htbp]
	\centerline{\includegraphics[width=0.9\textwidth]{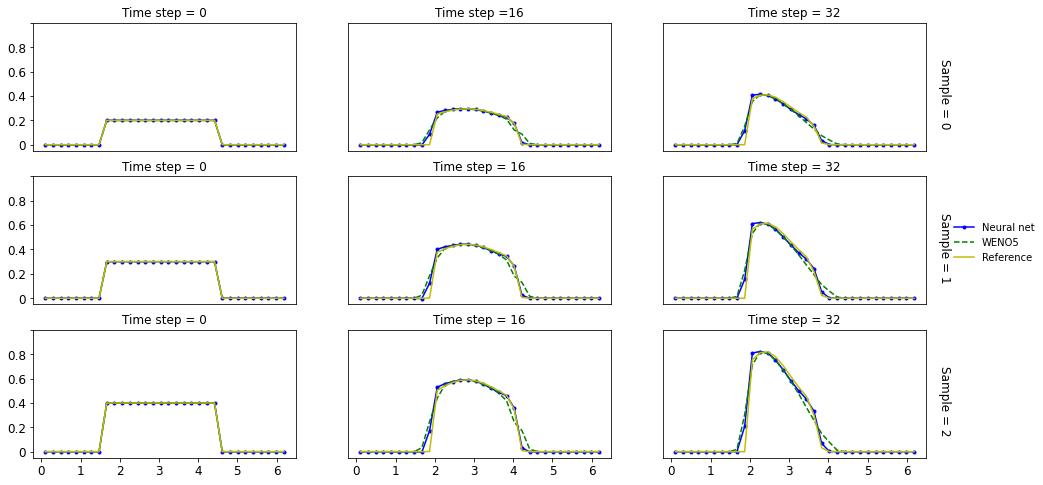}}
	\caption{Numerical solutions of three test samples for the transport equation with a variable coefficient in Example \ref{eg:variable}. CFL=0.6.}\label{fig:variable_three_curve}
\end{figure}

We first plot three test samples at several instances of time during forward integration with CFL=0.6 in Figure \ref{fig:variable_three_curve}. It can be observed that our solver can accurately resolve the solution structures with sharp shock transition and outperforms WENO5.

\begin{figure}[!htbp]
	\centering
	\subfigure[]{\includegraphics[width=0.45\textwidth]{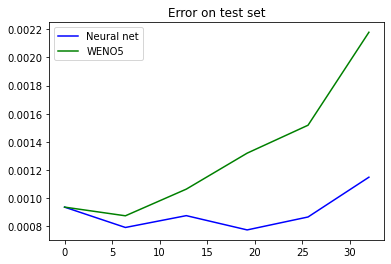}}
	\subfigure[]{\includegraphics[width=0.45\textwidth]{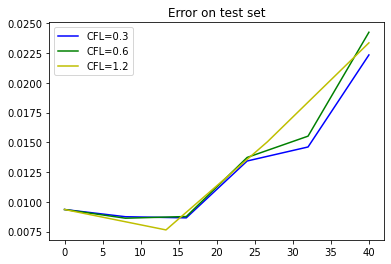}}
	\label{fig:variable_error}
	\caption{Time histories of the errors for the transport equation with a variable coefficient in Example \ref{eg:variable}. (a): Errors by the proposed method  compared with WENO5 averaged over all test samples. (b): Errors by the proposed method with different CFLs.}
\end{figure}

\begin{figure}[!htbp]
	\centerline{\includegraphics[width=0.48\textwidth]{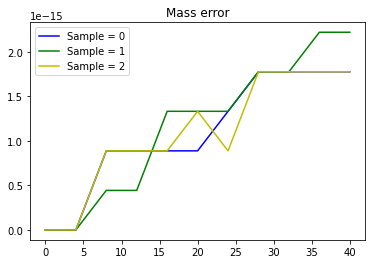}}
	\caption{Time histories of the deviation of  the total mass of three test samples for the transport equation with a variable coefficient in Example \ref{eg:variable}.}\label{fig:variable_mass}
\end{figure}

Figure \ref{fig:variable_error}(a) shows the time histories of the mean square errors, averaged over all test examples. The proposed method achieves the reduction of the error by a factor of 3.7 in comparison to the WENO5 method. We further report the time histories of errors for the method with three different CFL numbers in Figure \ref{fig:variable_error}(b), and it is observed the errors are comparable. As demonstrated in Figure \ref{fig:variable_mass}, the proposed method can conserve the total mass up to machine precision. 

\subsection{Two-dimensional transport equations}

In this subsection, we present the numerical results for simulating several 2D benchmark advection problem.

\begin{example}
	\label{eg:2dlinear}
	We solve the following constant-coefficient 2D transport equation 
	$$
	u_t + u_x + u_y= 0,\quad (x,y)\in[-1,1]^2,
	$$
	with periodic boundary conditions.
\end{example}

The training data is generated by coarsening 30 high-resolution solution trajectories over a
$256\times 256$-cell grid by a factor of 8 in each dimension, and each trajectory is initialized as a square wave with height randomly
sampled from $[0.5,1]$ and width from $[0.3,0.5]$. One trajectory contains 256 sequential time steps, and the CFL number is chosen within the range of $[0.3,1.8]$. We set the stencil size to be $5\times 5$. For testing, the initial conditions are sampled from the same range of width and height. The reference solution is generated by WENO5 with the mesh of $256\times 256$ cells and down sampled to the original coarse grid.

\begin{figure}[!htbp]
	\centerline{\includegraphics[width=0.9\textwidth]{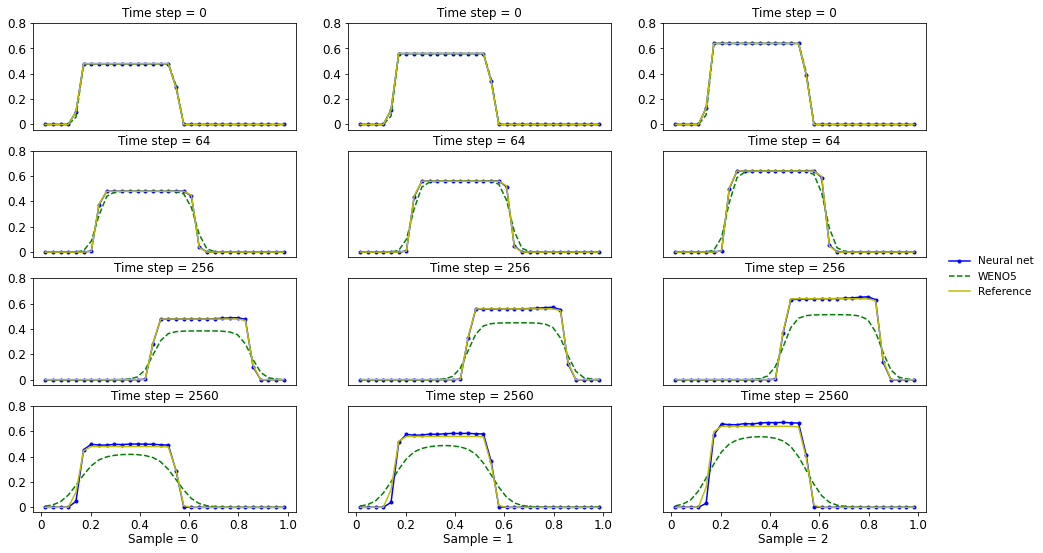}}
	\caption{1D cuts at $y=0.5$ of three test samples for 2D transport equation with constant coefficients in Example \ref{eg:2dlinear}.  CFL=0.6.}\label{fig:2d_three_curve}
\end{figure}

\begin{figure}[!htbp]
	\centerline{\includegraphics[width=0.9\textwidth]{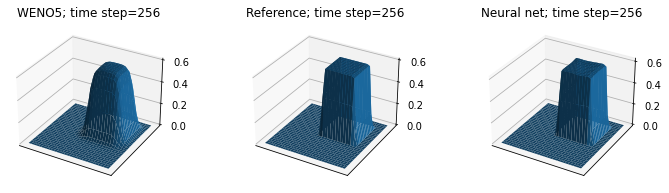}}
	\caption{One test samples for the 2D transport equation with constant coefficients in Example \ref{eg:2dlinear}. 2D plots of the solutions at time step 256. CFL=0.6.}\label{fig:2d_curve}
\end{figure}

Figure \ref{fig:2d_three_curve} shows 1D cuts of solutions at $y=0.5$ for three test examples at several instances of time during the forward integration with CFL=0.6. It is observed that the proposed method significantly outperforms WENO5 in resolving shocks sharply without introducing spurious oscillations. Furthermore, even after conducting simulations over a long period of 2560 time steps, the proposed method still produces highly accurate results. For a more effective comparison,  we also provide the 2D plots of the solutions at time step 256 in Figure \ref{fig:2d_curve}. It can be seen that the solution produced by WENO5 exhibits noticeable smeared shocks, whereas the solution by the proposed method exhibits much sharper shock resolution.


\begin{example}\label{sec:defor}
	In this example, we simulate the 2D deformational flow proposed in \cite{leveque_high-resolution_1996}, governed by the following 2D tranport equation
	\begin{equation}\label{eq:2d:defor}
		u_t + (a(x,y,t)u)_x+(b(x,y,t)u)_y = 0,\quad (x,y)\in[0,1]^2, 
	\end{equation}
	with the velocity field is a periodic swirling flow \begin{equation}\label{eq:deformational flow}
		\begin{aligned}
			a(x,y,t) &= \sin^2(\pi x)\sin(2\pi y)\cos(\pi t/T),\\
			b(x,y,t) &= -\sin^2(\pi y)\sin(2\pi x)\cos(\pi t/T).
		\end{aligned}
	\end{equation}
	It is  a widely recognized benchmark example for transport solvers. The solution profile is deformed over time.
	At $t=\frac{1}{2}T$ the direction of this flow reverses, while  the solution returns to the initial state at $t = T$, completing a full cycle of the evolution.
	
	
\end{example}

We set the period $T=2$ and the initial condition to be  a cosine bell centered at $[c_x,c_y]$:
\begin{equation}\label{eq:deformation-initial}
	\begin{aligned}
		u(x,y) &= \frac{1}{2}[1+\cos (\pi r)]\\
		r(x,y) &= \min\left[1,6\sqrt{(x-c_x)^2+(y-c_y)^2}\right].
	\end{aligned} 
\end{equation} 
We initialize 30 trajectories with $c_x$ and $c_y$ randomly sampled from $[0.25,0.75]$ using a high-resolution mesh of $256\times256$ cells, which are coarsened by a factor of 8 in each dimension as the training data. Each solution trajectory in the training data contains a sequence of time steps from $t=0$ to $t=T$.  We set the stencil size to be $5\times 5$. During testing, the initial conditions are sampled from the same distribution as the training data.

\begin{figure}[!htbp]
	\centerline{\includegraphics[width=0.9\textwidth]{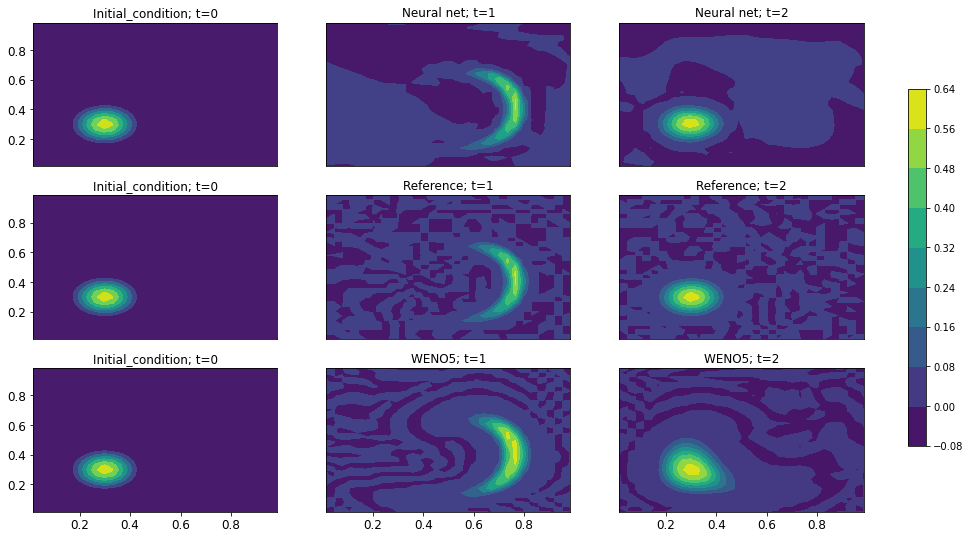}}
	\caption{Contour plots of the numerical solutions of the 2D deformational flow at $t=0,\,1,\,2$  in Example \ref{sec:defor} for one test sample. The solution profiles are significantly distorted at $t=1$, and then the flow reverses. At $t=2$ the solutions return to the initial states. }\label{fig:defor_onebell}
\end{figure} 

In Figure \ref{fig:defor_onebell}, we show the contour plots of the numerical solutions computed by the proposed method and WENO5 along with the reference solution for one test sample. Note  that the solution is significantly distorted at $t=T/2$ and returns to its initial state at $t=T$. It is observed that our method produces a result that is in good agreement with the reference solution. However, the solution obtained using WENO5 noticeably deviates from the reference solution due to a large amount of numerical diffusion.

\begin{figure}[!htbp]
	\centerline{\includegraphics[width=0.9\textwidth]{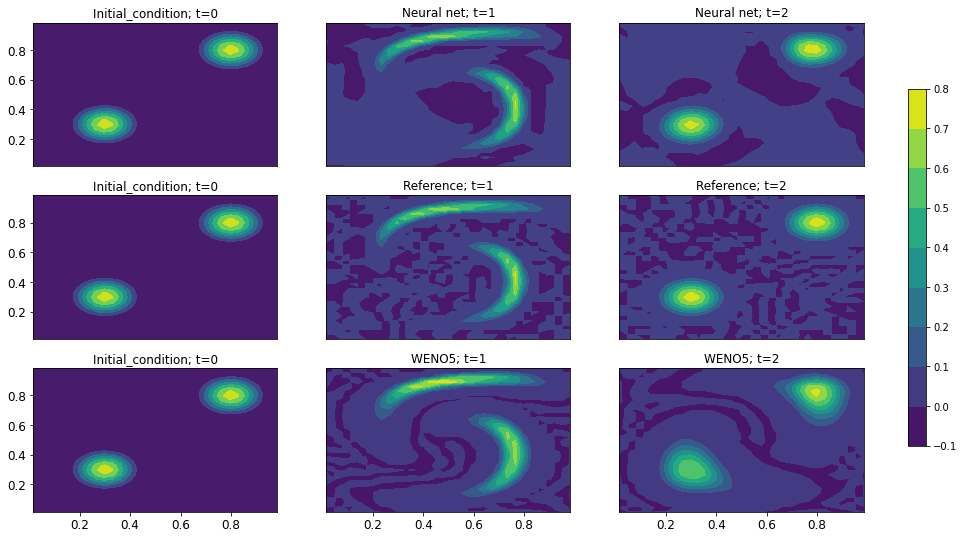}}
	\caption{Contour plots of the numerical solutions of the 2D deformational flow at $t=0,\,1,\,2$  in Example \ref{sec:defor} with two cosine bells. The proposed ML model is trained with a data set for which each trajectory only contains a single cosine bell.}\label{fig:defor_twobell}
\end{figure}

Furthermore, the model which was trained using data from solution trajectories featuring a single bell can  generalize to simulate problems with an initial condition containing two bells, as shown in Figure \ref{fig:defor_twobell}. The observation is similar to the single bell case. Last, we compare the errors of the numerical solutions by our method and by WENO5 at $t=T$ in  Table \ref{table:defor}. It is observed that the error of the ML-based SL FV method over the mesh of  $32\times32$ cells is  much smaller that of WENO5 with the same mesh size and is comparable to that of WENO5 over the finer mesh of $128\times128$ cells, demonstrating the efficiency of the proposed method.

\begin{table}[!htbp]
	\centering
	\caption{Accuracy comparison for the ML-based SL FV method and WENO5 of three test samples for the 2D deformational flow in Example \ref{sec:defor}.} \label{table:defor}
	\smallskip
	\begin{tabular}[c]{c c c c}
		\hline
		Samples \textbackslash Method & WENO5 ($128\times 128$) & WENO5 ($32\times 32$) & ML-based SL FV method \\
		\hline
		Sample = 0 &  9.269E-06 &  1.173E-03 & 1.955E-05 \\
		Sample = 1 &  1.069E-05 &  1.982E-03 & 9.753E-06 \\
		Sample = 2 &  9.302E-06 &  1.623E-03 & 1.107E-05 \\
		\hline
	\end{tabular}
\end{table}

\section{Conclusion}\label{sec:con}

In this paper, we proposed a machine-learning-assisted semi-Lagrangian (SL) finite volume (FV) scheme for efficient simulations of transport equations. Our method leverages a convolutional neural network to optimize SL discretization using high-resolution data, eliminating the need for costly upstream cell tracking. With a fixed 5-cell stencil, the CFL number can reach as large as 1.8. The inclusion of a constraint layer in the network ensures total mass conservation to machine precision. Numerical experiments show superior performance compared to the WENO methods. Future work includes extending the method to nonlinear transport equations such as the Vlasov system and investigating the use of graph neural networks for accommodating  unstructured meshes, adaptivity, complex geometries, among many others.  

\section*{Acknowledgments}
Research work of W. Guo is partially supported by the NSF grant NSF-DMS-2111383, Air Force Office of Scientific Research FA9550-18-1-0257. Research work of X. Zhong is partially supported by the NSFC Grant 12272347.

 \bibliographystyle{abbrv}

\end{document}